\documentclass{article}
\usepackage[a4paper, total={6in, 8in}]{geometry}

\usepackage{graphicx}        
\usepackage[bottom]{footmisc}

\usepackage{newtxtext}       %
\usepackage{newtxmath}       

\newtheorem{question}{Question}[section]
\makeindex             

\begin{document}

\title{Is math useful?}
\author{Alberto Saracco}
%
%
\maketitle

\abstract{\textit{Is math useful?} might sound as a trick question. And it is. Of course math is useful, we live in a data-filled world and every aspect of life is totally entwined with math applications, both trivial and subtle applications, of both basic and advanced math. But we need to ask once again that question, in order to truly understand what is math useful for and what being useful means. Moreover, is it knowledge of math useful for a class of specialists, or for political leaders or for all people at large? Being more on a concrete level, why does math need to have a central role in education?\\
Each section will be titled by a question. And each section will not give an answer, but ---at least I hope--- provide some food for tought to the reader, in order to try to come up with his or her own answers.\\
I feel that these kind of questions are at home in a book devoted to the interplays between mathematics and culture: what is the space we should give to math in culture and what is math's role in becoming a complete citizen?}

\begin{figure}[h]
		\includegraphics[width=0.60\textwidth]{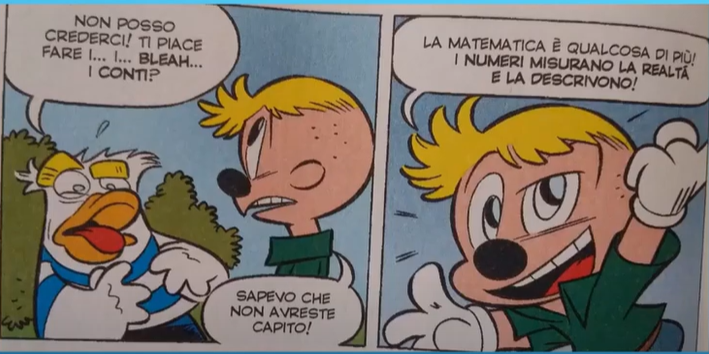}
	\caption{\textit{Math describes reality}, according to young Roby Vic, a Disney version of the ESA austronaut Roberto Vittori in the story \textit{3 gradini per le stelle --- Paperino Paperotto, Roby Vic e i conti... alla rovescia} (\textit{3 steps to the stars --- Donald Duckling, Roby Vic and the countdowns}) \cite{RV}}
	\label{figRV}
\end{figure}

\section{Is math useful?}
\begin{quotation}  The mass of mathematical
truth is obvious and imposing; its practical applications, the
bridges and steam-engines and dynamos, obtrude themselves on
the dullest imagination. The public does not need to be convinced
that there is something in mathematics.
All this is in its way very comforting to mathematicians, but it
is hardly possible for a genuine mathematician to be content with
it.

\begin{flushright}A mathematician's Apology \cite{H} \S 2 --- G. H. Hardy
\end{flushright}
\end{quotation}

\begin{quotation} [...] the most ‘useful’ subjects are quite
commonly just those which it is most useless for most of us to
learn. It is useful to have an adequate supply of physiologists and
engineers; but physiology and engineering are not useful studies
for ordinary men[...]

\begin{flushright}A mathematician's Apology \cite{H} \S 20 --- G. H. Hardy
\end{flushright}
\end{quotation}

The English mathematician Hardy dealt very well with the subject of this paper, and I will often cite his famous \textit{Apology}, written over 70 years ago. Since we are dwarves sitting on the shoulders of giants, I hope I'll be able to see a little further and give some new ideas on the subject.

More precisely I would like to deal with the problem posed by the above two quotes: no one is usually fool enough to deny the usefulness of mathematics to our society, but the usefulness to a society is not at all the same as the usefulness to an individual. \textit{What is math useful to me?} will be our next to final section.

Before getting there, tough, we have a long way. We first have to understand what is the usefulness of math and how math is (and can be) used.

\section{How is math used in war time?}
\label{sec:WWII}

 \begin{quotation}Ten, twenty, thirty, forty, fifty or more

The bloody Red Baron was running up the score

Eighty men died trying to end that spree

Of the bloody Red Baron of Germany

\begin{flushright}Snoopy vs the Red Baron (1966) ---  The Royal Guardsmen\footnote{Usually the books of the series \textit{Imagine Math} are the proceedings of the meeting on mathematics and culture held in Venice. This year, due to the pandemic of Sars-Cov2, the meeting has not taken place. Anyhow you may think of me beginning my lecture playing with planes while the song \textit{Snoopy vs the Red Baron} is being played. I suggest you to listen to this song while beginning to read this chapter, to put yourself in the right mood.}
\end{flushright}
\end{quotation}

\begin{figure}[h]
		\includegraphics[width=0.40\textwidth]{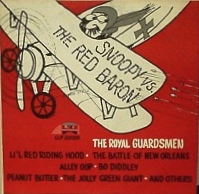}
	\caption{The cover of the album Snoopy vs the Red Baron (imagine from wikipedia)}
	\label{figRedBaron}
\end{figure}

Math has always been considered a strong ally in war time. Archimedes used math (and physics) to construct parabolic glasses in order to set on fire Roman ships. Math has been used to compute the line of firing of cannons: modern ballistics was born due to an English mathematician, Robins, who in 1942 wrote New Principles of Gunnery, a treaty which was used till World War II.  Mathematicians have always been considered precious for war and enrolled for their logical and computing abilities. The English mathematician Littlewood, closed friend and fellow mathematician of the already cited Hardy, served in the Royal Garrison Artillery during World War I.

Up to World War I, tough, the math used in war time was quite elementary: basic geoemetry and physics.

During World War II math played a fundamental role and many different areas of mathematics turned out to be useful for winning the war. Everyone knows the story of how Alan Turing cracked Enigma, the Nazi cryptography tool, and heavily contributed to lead the Allies to a victory. Both cryptography and decryptography are based on deep math.

Moreover there is some statistics in figuring out the number of tanks produced by the Nazi: the problem of estimating the maximum of  a discrete uniform distribution from sampling without replacement. In simple terms, suppose there exists an unknown number of items which are sequentially numbered from 1 to N. A random sample of these items is taken and their sequence numbers observed; the problem is to estimate N from these observed numbers. This problem is called the German tank problem, since it was of uttermost importance to the Allies: they wanted to estimate the number of German tanks just by knowing the serial numbers of the few tanks captured.

\begin{figure}[h]
		\includegraphics[width=0.40\textwidth]{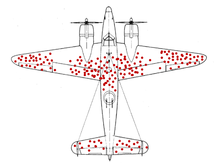}
	\caption{Area of damage of damaged airplanes returned to base during WWII (imagine from wikipedia).\\
Where do the airplanes need to be reinforced?}
	\label{figBias}
\end{figure}

But probably the nicest use of math in WWII is that Wald did to evaluate where planes would need additional armor against enemy's shootings. Data of damaged airplanes was collected by US military, leading to the picture of figure \ref{figBias}. The US military concluded that the area in need of ticker armor where the ones with the most shoots. Wald concluded the opposite: the areas in need for additional armor were precisely the one with least red dots, since the sample was made up just of planes who survived the enemy's fire. The planes which were shot elsewhere did not made their way back home: they simply crushed down, as if they were shot by the bloody Red Baron. So the parts to reinforce were exactly those which when shot did not allow the plane to come back and its damage to enter the stats.

Math is no doubt useful in war time, both for computational purposes (and in this I put also physics and computer science) and for the mathematical-logical thinking.

\begin{question}{What is math useful for?} I have no doubt that the reader, when reading the title of this chapter, immeditely tought that math \textit{is useful} and did not think about math in war time, which ---depending on how you put it--- may be described as useful or bloody dangerous. So, why did I choose such a subject to begin my chapter? I'll let Hardy answer:

\begin{quotation} I once said that ‘a science is said to be
useful if its development tends to accentuate the existing inequalities in the distribution of
wealth, or more directly promotes the destruction of human life’

\begin{flushright}A mathematician's Apology \cite{H} 21 --- G. H. Hardy
\end{flushright}
\end{quotation}

This extremely pessimistic phrase was spoken by Hardy in 1915, when times were dark and there was little space for hope and for the future of humankind. Nevertheless, way too often the \textit{usefulness} of something has indeed had the effect to accentuate inequalities or favour wars, as Hardy stated. Hence, I feel that we should start discussing usefulness of math by starting from its darkest sides, not hiding them under the carpet, but being well conscious of their existence.

\end{question}

It is also worth noticing that, usually, when a war time example of usefullness of math is made, it is usually a situation in which the good Allies used math to win against the Nazi.

It is as if, when talking about the applications of math to the real world, we try to hide the darker sides of math's applications, and ---if ever we talk about math and war--- cite only occasion in which math has been used to make the Allies won WWII vs the Nazi, i.e. show war-time-math as the hero in a classic war movie.

Of course, reality is much more various than a movie, and in war math helped killing people as much as saving them.

We must deal with this whenever we try to answer the question whether math is useful or not: math is a tool, and like most tools it can be used in a wise or a wicked way.

\section{How can pure math be useful?}
\label{sec:2}

\begin{figure}[h]
\center
		\includegraphics[width=0.75\textwidth]{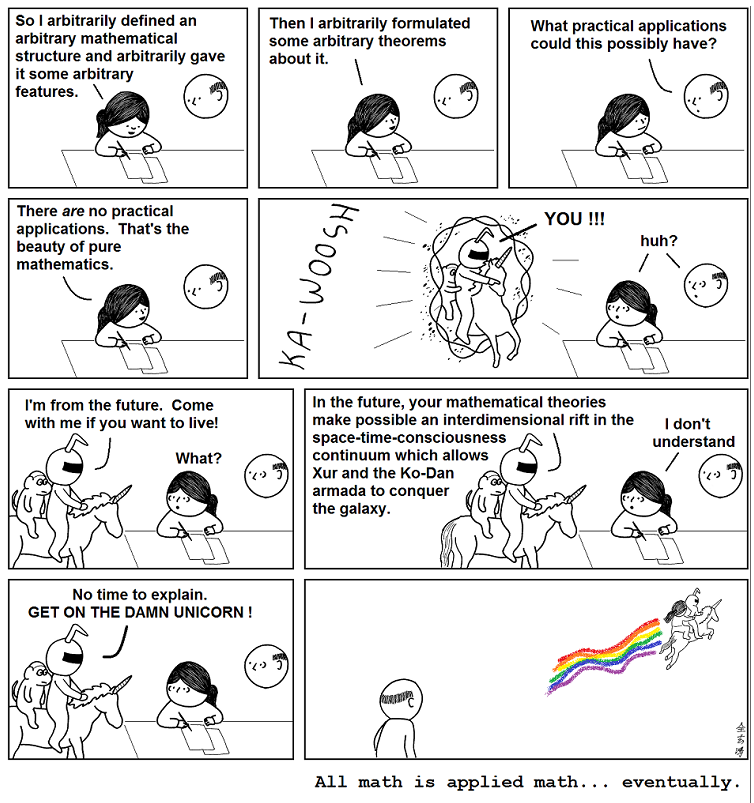}
	\caption{All math is applied math... eventually \cite{AG}}
	\label{figAG}
\end{figure}

There are uncountable\footnote{Obviously this is an hyperbole, since everything in the real world is not only countable, but finite!} applications of math in everyday's life. While most of the applications known to the wide public rely on basic math, or on math born explicitely for applications, I would like to give some example of pure math which later on turned into applied math. As the online comics Abstruse Goose puts it, \textit{all math is eventually applied math} (see figure \ref{figAG}). 

If you feel a geek online comics is not a good enough reference\footnote{And you'd be wrong. Comics are totally part of cuture and they also had a good place in this series of books, see e.g.\ \cite{A,S}}, I'll go with Galileo:

\begin{quotation}
[Nature's book] is written in mathematical language, and its characters are triangles, circles and other geometrical figures.\footnote{\textit{[Il libro della natura] è scritto in lingua matematica, e i caratteri son triangoli, cerchi, ed altre figure geometriche}, in the original}

\begin{flushright}\textit{Il Saggiatore} --- Galileo Galilei\end{flushright}
\end{quotation}

Our limits in applying mathematics in describing the world are just those of our knowledge (deep and true knowledge) of it\footnote{Is it "the world" or "math"? I'll leave the answer to the reader}. We may think our knowledge as a box of tools. As soon as we have a tool in it, we may find some uses for the tool. When we do not have the tool (or we even ignore its existence), we cannot find uses for the tool. And mathematical tools, being completely general and abstract, have a wide range of possible applications. What is really needed is to create the mathematical tools (doing mathematical research, both pure and applied) and handle them to people who may need them.

Thus, in this section I go with some example of pure math applications to the real world.

\subsection{Number Theory and Cryptography}
\label{subsec:2}
The Abstruse Goose comics (figure \ref{figAG}) gets it right: no matter how pure and far from applications a part of mathematics may be, once it is in our tool box it is only a matter of time since it will find some applications\footnote{Again an hyperbole: some theorems are just too weird to find an application... or are they just too weird for now? Maybe just because we do not understand that result well enough?}.

\subsubsection{Number Theory} So it is just appropriate to begin with an example about number theory, the field of research of Hardy, who was absolutely proud of doing research in a pure area of math, with no applications whatsoever.

Hardy in his \textit{A mathematician's apology} \cite{H} makes quite a point of personal pride in number theory being a completely pure (and useless) branch of mathematics. And he was quite right! Number theory deals with the distribution and the properties of prime numbers and so it is a subject of great charm, ancient (the proof that there are infinite prime numbers and Eratosthenes' sieve date back to Ancient Greece) and full of elementary problems (e.g. Goldbach's conjecture), which are easy to state and extremely difficult to solve. This characteristic lead to a heap of amateurs trying to solve very difficult conjectures in this field across centuries\footnote{And this often lead to frustration professional number theorist who continue to recieve "proofs" from amateurs, see e.g the nice \textit{Dialogue on Prime Numbers} written by Zaccagnini \cite{Z}}. Few succeded, most didn't, and a lot of appearently simple conjectures are still unsolved.

So Number Theory always had a great appeal, both to professional mathematicians, amateurs and the wide audience. But no one ever questioned its being a totally pure and abstract area of mathematics, whose interest belonged all to the world of pure ideas and not to our material world.

At least, that was so until computer age begun and some old Number Theory theorems by Fermat become useful for cryptography.

\subsubsection{Cryptography} 
\begin{figure}[h]
\center
		\includegraphics[width=0.50\textwidth]{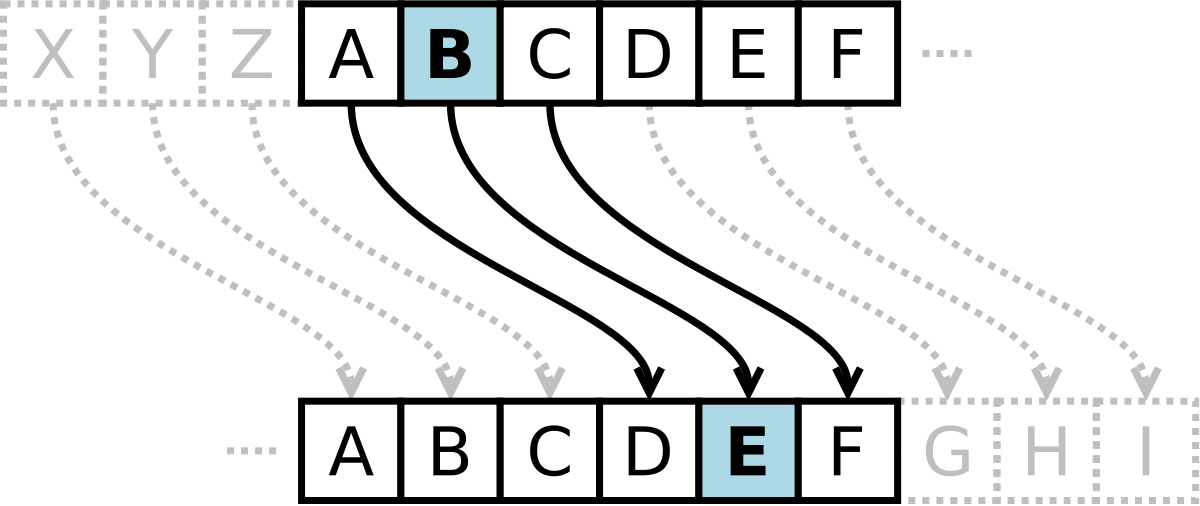}
	\caption{Caesar's cryptographic method: just replace a letter with the one three places after in alphabetical order (image from wikipedia)}
	\label{figCC}
\end{figure}

Also Cryptography is an ancient subject. Sending secret messages has always been of crucial importance in war time (as we already stated in section \ref{sec:WWII}) and the first use of Cryptography dates back to Julius Caesar, who sent messages replacing each letter with the one 3 places after that in alphabetical order (see figure \ref{figCC}).

To read the original message, one should just reverse the arrows and replace each letter with the one 3 places before in alphabetical order.

This method of crypography has several problems, of course.

\begin{figure}[h]
\center
		\includegraphics[width=0.60\textwidth]{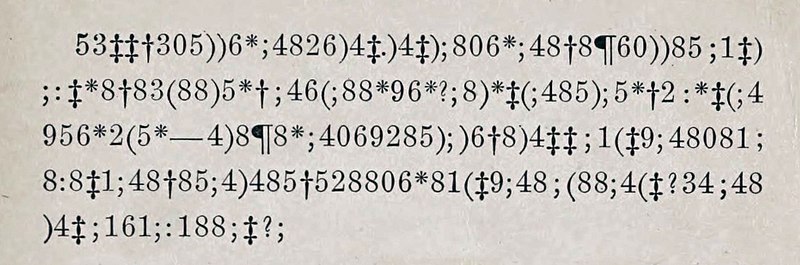}
	\caption{The cryptogram in Poe's novel \textit{The Gold-Bug} \cite{GB}.}
	\label{figGB}
\end{figure}

First of all, if one knows how to encrypt a message, also knows how to decrypt the message. Secondly, the possible shifts are just one less then the number of characters in your alphabeth (not counting the 0-shift which does not encrypt): not many to check fast by hand. Moreover, even if not a simple shift is used to encrypt but rather any permutation, if the message is long enough, a simple statistical analysis of most frequently occurring characters may yield to an easy decryptation (as it is done in E. A. Poe's \textit{The Gold-Bug}, see figure \ref{figGB}). Finally, both the receiver and the sender must know the crypting and decrypting keys in order to have a crypted communication between them. And how do they exchange the keys?

With the computer era, decrypting messages has become easier and easier. The faster the computers, the better encryption methods had to be.

The breaking of the Nazi encrypting machine \textit{Enigma} by a huge group lead by Alan Turing was a key turning point of WWII.

\subsubsection{Number theory and cryptography}
Number theory was used to solve the last and biggest problem in cryptography. Namely, number theory allows for a method in which the encryption key is public but the decryption key is private, thus allowing anyone to send messages to the receiver (e.g. your password to a website or the PIN of your credit card to your bank) without having to worry about a third party decrypting the message.

This has been a mayor breakthrough in cryptography and its applications are amusing.

From the theoretical point of view, the RSA method is really simple. One needs to find two distinct prime numbers, $p$ and $q$ and computes $n=pq$ and $m=(p-1)(q-1)$. Then a number $a$ such that $(a,m)=1$ is choosen and the number $b$ such that $ab\equiv 1$ (mod $m$) is computed.

The couple $(n,a)$ is the public key and is known to everyone. The number $b$ is the private key and it is secret. The message is translated into a number $x<n$ and the sender sends the encrypted message $y<n$, where $y\equiv x^a$ (mod $n$). The receiver computes
$$y^b \equiv (x^a)^b = x^{ab} \equiv x \  (mod\ n)$$
thanks to Euler Theorem, thus getting to know the original message.

The operations of taking a power up to a congruence class is not much time consuming and can be easily done by a computer. Finding out from $n$ its prime factors $p,q$ is completely a different task, in term of computation time. Of course, the greater the computational power of computers, the bigger the two primes $p$ and $q$ need to be. Nowadays the RSA key is 128 bits long (or 256 bits for TopSecret tasks).

\subsection{Radon-Nikodim antitransformation and computed axial tomography}

Looking inside a body may be a difficult task. Our body is not transparent and cutting a person in order to see what's on the inside may not always be a good idea. Radiography, using X-rays, helped in seeing bones, since the rest of the body is trasparent to X-rays, but they are not a good means to inspect soft parts of our bodies.

Medicine was in search for a tool we were appearently lacking: a way to see inside our bodies without tearing them apart. The tool was only appearently lacking. Indeed math has invented ways to transform local informations into global ones and viceversa: trasformation and antitransformation. There are several of them, and they answer to different needs, but actually what a transformation does is taking as an input a function or a series of numbers, and giving back another function or a series of numbers; the antitrasformation going the opposite direction and being an inverse to the transformation. Usually these tools work computing integrals.

Sending rays through your body and see how much they were absorbed was not a new idea (indeed it was used with X-rays and radiography), but it is just in the early Seventies that a physicist (Allan Cormack) and an engineer (Godfrey Housefield) had the idea of using Radon-Nikodim transformation and its inverse in order to compute from the information of rays absorbed in the various direction a 3D model of the inside of a body. This application of math eventually led to the Nobel Prize for Medicine (in 1979) for the two and gave a huge tool of diagnosis to hospitals all around the world.

When Radon-Nicodim transformation was developed in 1917, it was a completely pure and "useless" tool of high mathematics. Of course, computers were far from being invented, at that time, and practical applications of the Radon antitransformation were unforeseeable. But, as we said, all mathematics is eventually applied mathematics. And you never know when something you know may actually turn useful. In any case, it is better to know more than less.

\section{Why politicians should know math?}
All of the above is a bunch of examples showing what is math useful to us as a species, as a community. But of course, we may ask ourselves if math knowledge should not be simply limited to mathematicians, engineers abd other people who may use it in their work for the benefit of the community at large. After all, we do not need people know exactly how a bridge is constructed, how a TV works, or how to repair a broken engine. For that, we use people who know how to do. Why should math be different?

I will address this question in the following sections. First, let us consider why politicians should know math.

Our modern world is a world filled with data and numbers. Decisions must be made based on those numbers and those data. But interpreting data is far from obvious, as the "survivor's bias" example should have already clearified. An inability to correctly interpret data may turn into a disaster. Indeed statistics is quite difficult.

\subsection{Education system in the US, covid-19 death toll and Simpson's paradox}
\begin{figure}[h]
\center
		\includegraphics[width=0.60\textwidth]{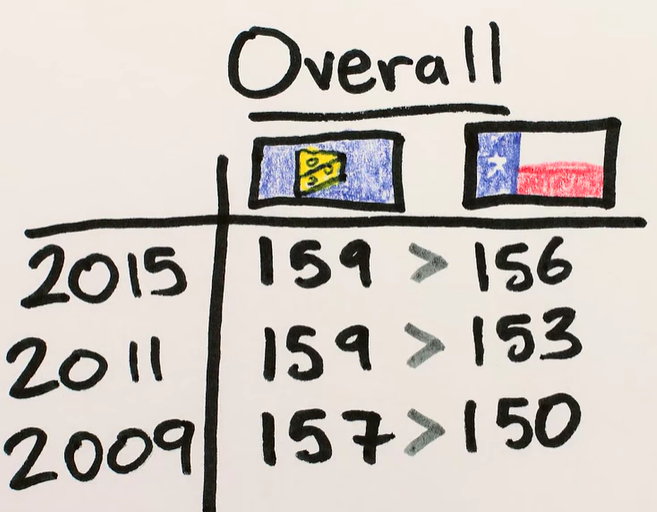}
	\caption{Data of Texas and Winsconsin overall results in standardized tests  (source minutephysics \cite{MP}).}
	\label{overall}
\end{figure}

For many years, Winsconsin's students performed consistely better than Texas' students in standardized tests (see figure \ref{overall}).

One could conclude that Wnsconsin's education system is way better than Texas', and a politician willing to improve Texas' education system may be tempted to copy the one of Winsconsin. But is that a good idea?

Knowing the mean performance of a huge number of students for a long time may sound as pretty solid evidence towards this claim. But statistics is full of surprises.

\begin{figure}[h]
\center
		\includegraphics[width=0.60\textwidth]{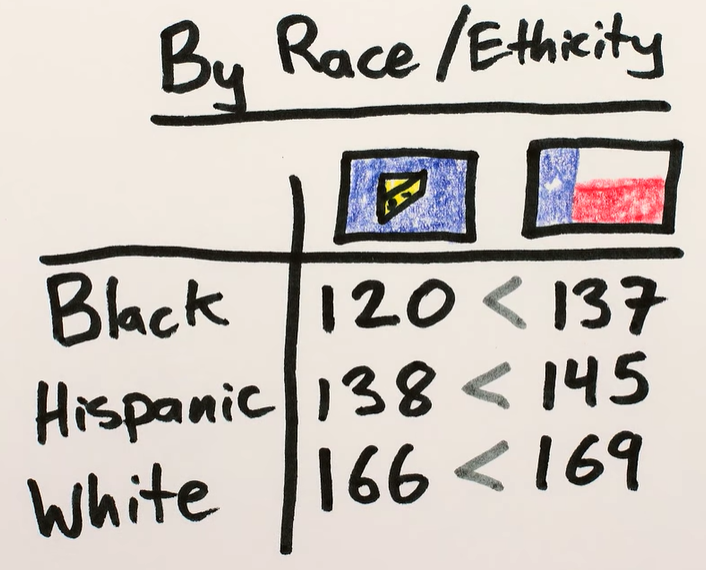}
	\caption{Data of Texas and Winsconsin results in standardized tests, divided by race (source minutephysics \cite{MP}).}
	\label{race}
\end{figure}

Namely, if we divide the data of the students of the two States among different ethnic group (and we know that ethnic group correlates with wealth which correlates with education level), a surprise pops out: white Texas students outperform white Winsconsin students, black Texas students outperform black Winsconsin students and hispanic Texas students outperform hispanic Wisconsin students (see figure \ref{race}). 

So it actually looks that, when seen broken by race, data suggest that Texas' education system is better than Winsconsin's. How can data tell two different things? First of all, one of the problem is that the mean of some data is not the same as the mean of the means: it depends on how many data are there in every subgroup in which data have been divided. Wisconsin's population is much whiter than Texas': thus the overall mean of Wisconsin is much more tilted towards the white mean (which is the etnic group performing better in the test) than it is the mean of Texas.

\begin{figure}[h]
\center
		\includegraphics[width=0.60\textwidth]{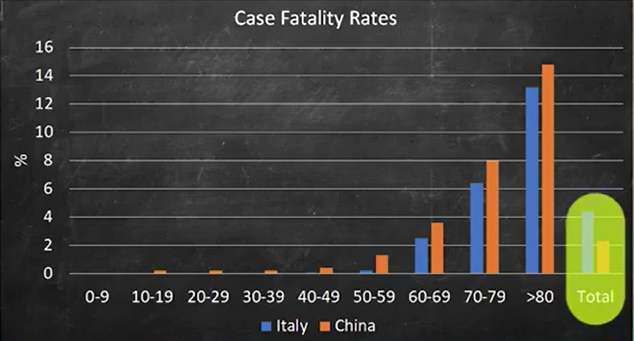}
	\caption{Data of Italy and China fatality rates of Covid-19 cases (image from \cite{B}).}
	\label{CovidFatality}
\end{figure}

A similar situation happened when comparing the death toll in Italy and China at the beginning of the Covid-19 pandemic. Indeed, the overall fatality rate of the disease in Italy was bigger than the overall fatality rate of the virus in China, but ---when people infected with Sars-CoV2 were split in age groups--- in every single age groups the fatality rate was greater in China than in Italy (see figure \ref{CovidFatality}).

In this case the problem is that in Italy there were much more old people sick with Covid-19 than there were in China. And Covid-19 has a higher fatality rate in older patients. So, this explains the appaerent discrepancy of the data.

This phenomenon, where there is a positive correlation overall, while ---when data is divided in groups--- there is a negative correlation, is called Simpson's paradox. Simpson's paradox is one of the thing a politician should be aware of, before taking action according to data.

But actually the problem is deeper than that. Indeed, is Wisconsin's school system better or worse than Texas'? Is this second way to see data the correct one? If we say that etnicity correlates with wealth (or with parents' education level) and this last thing correlates with results in tests, why are we using etnic group and not wealth (or parents' education level) directly in order to interpret our data? A third look at data may be needed.

The point is: if you believe data are objective and need not to be interpreted and anylized with a close look, you are likely to be fooled by data. If you know how statistics (and math) works, you are more likely not to get fooled and to take a second (or even a third) look at data before taking action (and possibly going in the wrong direction).

This is why politicians should have a good base in math and know how to analyze data. Before making decisions, the least you must have is correct data and infos, and possibly understand correctly what they mean.

\section{How do politicians use math?}
You might say that politicians do not really have to know and understand math, in order not to fall into such errors, but just to have good advisors who do know math. And indeed they have. Plenty of them. And here we get to the problem.

First, as we have seen, data and numbers are far from being \textit{objective}: data must be interpreted and investigated in order to understand what they say, but they are also easily bent to furnish support to almost any political view. So, sadly, way too often the scientific advisors of politicians try to cherry pick data or to present data in such a way to give a scientific-looking aspect to the political ideas they want to communicate. This when data are not right-away invented. But cheating too much is not even needed: the same data, presented with different words, from a different point of view, can lead to very different conclusions. And we must bear in mind that politicians often have a very skilled adving group whose only purpose is to find the best way to present data.

Another new interesting tool of math (or computer science) often used by politicians is given by big data: \textit{sentiment analysis} and \textit{trending topics} are fundamental in political communication. In our modern world we have an incredible amount of data about almost everyone: use of credit cards, posts or comments on social networks, our GPS position in real time, the shopping habits (both on-line and in physical stores), internet usage... The math of big data can extract patterns out of all this huge amount of data. And this is how your phone can suggest you the fastest path to go back home or where to buy the book you really want to read or the item you really needed. This can be useful, but of course all this information can be (and is!) used to make enormous profits.

Politicians are informed real-time about the hottest arguments of discussion (\textit{trending topics}) and on what most people think about the argument (\textit{sentiment analysis}), and so they are ready to band-wagon on the hottest topic with the coolest opinion. It almost does not matter whether the opinions expressed are coherent with one another or not: what is really important is to say something on the trending topic of the day, with their opinion being shared and viewed by the highest possible number of people. In time of election, people will recognise your name, and you'll have bigger chances of being voted, hence more votes. This is the core of marketing, applied to politics. Not so great if you think politics should be about solving the community's problems, but that's how it is. And there is a lot of math in that.

\subsection{Paradoxes of elections}
So, we have decided that politicians' biggest task is getting elected (even if they are really interested in doing their work for the benefit of the community: in order to do that, they need to be elected). Alas, the outcome of elections is far from being determined from what voters think, and the electoral system is crucial for the result. This is exactly the reason why politicans spend so much time discussing the electoral system. This subsection is mainly based on my paper \cite{SS}, on mathematical paradoxes of elections. I refer the reader to that paper  for greater details.

\begin{figure}[h]
\center
		\includegraphics[width=0.6\textwidth]{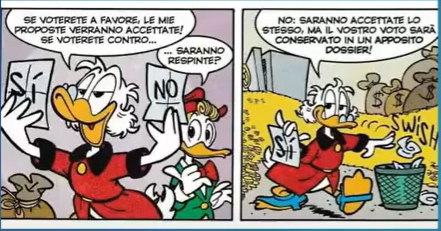}
	\caption{The only election satisfying Arrow's axioms is a dictatorship:\\
\textit{US: "If you vote YES, my proposals will be accepted. If you vote NO..."\\
GG: "they will be rejected?"
US: "No, they will be accepted, and your vote will be kept in an appropriate dossier"} \cite{NG}.}
	\label{Paperone}
\end{figure}

Unluckily, no electoral system is perfect. In 1951, the economist Kenneth Arrow \cite{KA} considered a very general definition of electoral system as a function (which he calls \textit{social choice function}) from the individual preferences among the alternatives of the electors to a single preference of the social group, where a \textit{preference} is a total ordering of the alternatives. Arrow introduces three desiderable properties of the social choice function:
\begin{itemize}
\item[\textbf{A1}] \ (sovranity of electors) The function is surjective, i.e.\ if the electors agree on the desired outcome, they can vote (choose their individual preferences) in order to have that outcome;
\item[\textbf{A2}] \ (positive correlation) If in a certain situation the social choice function says $x$ is better than $y$, in any other situation in which the only change in any elector's preference is that their ranking of $x$ gets higher, then $x$ is still better than $y$;
\item[\textbf{A3}] \ (Invariance under irrelevant alternatives) The relative position of $x$ and $y$ according to the social choice function depends only on the relative positions of $x$ and $y$ for each elector and not on the opinion on a third alternative $z$.
\end{itemize}
Arrow then proves that if there are at least three alternative, the only social choice functions satisfying the above three axioms are dictatorships: the social choice function is simply a projection on one on the factors, or differently said the "will" of the people is the "will" of a single individual, the dictator (see figure \ref{Paperone}).

This theorem actually means that, if more than 3 alternatives arre allowed, our electoral system must not satisfy one of the above properties if we do not want to have a dictator. Electoral systems usually do not satisfy the axiom \textbf{A3} and the outcome of election can be modified by the presence or absence of otherwise totally irrelevant political forces.

Politicians (or their advisor) well know this fact, and this explains all the fight on whether some little meaningless party should be allowed or not to partecipate in an election.

The Theorem of Arrow is based on Condorcet paradox, i.e.\ situations where you have 3 alternatives A-B-C and A wins vs B, B wins vs C and C wins vs A, in a rock-paper-scissor way.

The Theorem of Arrow may suggest that a system with only two alternatives to chose from is the best one. The most useful electoral system to force politicians to gather in only two major party, thus having a system with only two alternatives and a way out of Arrow's paradox, is a one-discrict one-seat system, where the party who gets the most vote in the discrit takes the seat, and all other votes for the other parties are meaningless.

Unluckily the system one-district one-seat has one big weakness: the outcome of the election strongly depends on the shape of the districts and a party which has the power to decide the shape of the districts may win the election in a 1 vs 1 race with as little as slightly more than $25\%$ of the votes. This is due to the fact that the party can lose $0$ to $100\%$ in slightly less then half the districts and win by barely one vote over $50\%$ in the remaining (slightly more then half) districts, thus winning the election. Of course it is impossible to have complete information about votes, but big data analysis gives parties a quite good level of knowledge about voting intentions, thus allowing an easy win even with if the electorate strongly favours the other party.

\begin{figure}[h]
\center
		\includegraphics[width=0.60\textwidth]{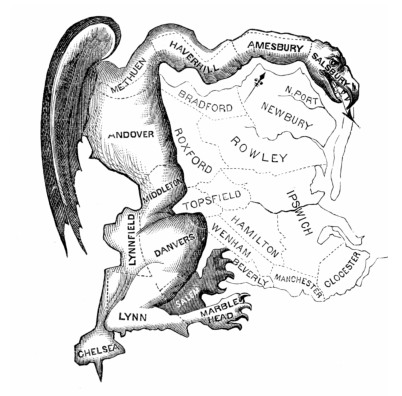}
	\caption{The satirical panel, with the salamander-shaped district, published on the Boston Centinel to mock Governor Gerry in 1812.\\
 Image of public domain, from wikipedia.}
\label{Gerry}
	\end{figure}

This art of carefully shaping the districts in order to win is called \textit{gerrymandering} in honour to the salamander-shaped district designed by Governor Elbridge Gerry to win an election (see figure \ref{Gerry}), but is still well used nowadays in the US, by both parties and districts of really weird shapes are not at all uncommon. Mathematical research on the subject of gerrymandering is very active, to limit gerrymandering, both by finding the objective subdivision on districts or by giving measures to find out whether there has been some gerrymandering going on in order to cheat or if the subdivision is fair.The two main approaches to the problem are an analytical approach using isoperimetric-like techniques and a discrete geometry approach, using weighted graphs. I proposed an approach of this second kind \cite{SS19}.

All this is why it is usually forbidden to change the rules of an election (or the shape of discrits) too near to the upcoming election. But of course, regulations do not completely stop politicians to use math for their own benefit.

\section{Is math useful to me?}

Let us now address the main point of this paper: how is math useful to \textit{me}? Why should \textit{I} learn math? Can't just a few people be knowledgeable of math for the whole society's benefit?

I once heard that when a kid learns to read by itself, it does no longer have to depend on others to read and can find out what's written around without having to trust others to read correctly to them. Math is a powerful tool to read the world, and knowing math gives you the power to independently analize the complexity of the world, without blindly trusting others to do that for you.

Be aware! I am not saying you should not trust others or the scientific community, not at all! What I am saying is that, being able to do the math by yourself or ---even better--- to mathematize a problem and to take a look at it through the magnifying glass of math, is always a good idea, not only to find out a correct answer yourself, but mainly to find out who is trying to trick you and who you should trust.

So a first answer is that the more you know, the less you are likely to be fooled by people who want to gain something by fooling you, by politicians or lobby who want to push their own ideas, or simply by arguments which may look plausible until inspected closely. Knowing math, or even better being able to reason with a math-oriented mind, is fundamental for every single individual.

Of course the problem is that it is not enough being able to reason correctly and foresee what's coming, if the majority of the population does not share this ability and is easily fooled into non rational behaviours. This is always true, but even more in a period where acting fast and correctly is the key to avoiding a disaster.

In October 2020, at a EU meeting, the German Cancellor Angela Merkel said "\textit{Once we got to this point, closures are the only possible choice: we should have acted earlier, but people would have hardly understood. They need to see hospitals' beds full...}"\footnote{Mine translation}. Angela Markel has a Ph.D. in physical chemestry and knew pretty well what was going on. She even gave, some month before, a very nice explanation of the meaning of the index $R_t$ in an epidemics. But that was not enough. She was knowledgeable, she was powerful, but still she could not act without her citizens being fully aware of what was going on. She was in the sad situation where she could foresee hospitals' beds getting full and death tolls raising to very high levels without being able to act tempestively.

The consequences of a low mathematical literacy of the vast majority of the population are terrible: more deaths, more pressure on the health system, more economical consequences... In order to avoid this in the future (the damages of the present situation cannot be undone, alas), a wide-spread math literacy is needed. Math allows to see you that an exponential growth of an epidemics means the hospitals' beds will be full and to act timely, in order not to have them full.

Without people's clearly understanding that, the action needed to solve the problem is also the action that will made people shout "\textit{Nothing happended! There was nothing to be worried about! We should have not done this}".

So, not only math gives you instruments to understand, analize the world and not getting fooled, but a wide-spread knowledge of mathematics will turn into a huge benefit for the whole society.

\section{Is math usefulness relevant to learners?}
I hope we have cleared that math is useful to everyone and to the population at large. Given that, the fact that something is useful to someone, it does not straightforwardly imply that it will be interested in learning that, much more when you are dealing with children or kids or young adults. Showing the utility of something may make want some students to learn something, but most of them will be simply bored as hell.

I will just quote Paul Lockhart, who makes his Salviati go right to the point\footnote{Bolding mine}

\begin{quotation} It may be true that you have to be able to read in order to fill out
forms at the DMV, but that’s not why we teach children to read. We
teach them to read for the higher purpose of allowing them access to
beautiful and meaningful ideas. Not only would it be cruel to teach
reading in such a way— to force third graders to fill out purchase
orders and tax forms— it wouldn’t work! We \textbf{learn} things \textbf{because
they interest us now, not because they might be useful later}.
\begin{flushright}A Mathematician's Lament \cite{L}  --- P. Lockhart
\end{flushright}
\end{quotation}

We learn because we are interested, because we are amused by something, not because we have to or because it will do us some good or it will make us a better person!

Luckily math is filled with interesting ideas and theories! We must never forget this, and when trying to appeal a young learner \textit{usefulness} should not be our guide through mathematics. Math was not born \textit{because} it was useful, but beacuse it is amusing. Math is filled with interesting problems, which can be given to kids and adult of different ages and knowledge, in order to hook them into mathematics.

For the sake of completeness, I will just present an example, but many more can be found in the very indicated reading of Lockhart's pamphlet \cite{L}.\vspace{0.2cm}

\textbf{Second degree polynomials.} Usually, when studying second degree polynomials, students are presented with a huge nomenclature (pure polynomial, spurious polynomial...) and a vast casistic to solve particular polynomial equations, then they are given a lot of exercises to practise the method that was given. After that, they are given the general formula
$$ax^2+bx+c=0 \ \ \Leftrightarrow\ \ x=\frac{-b\pm\sqrt{b^2-4ac}}{2a}$$
(possibly also with a second variant in case $b$ happens to be even) and then a new round of dumb exercises, each one equal to the previous. Totally boring.

I mean, I know that during the history of math all these different kinds of equations were solved (and given funny names), but math is not about zoology or funny names: math is about the struggle to find a path that lead to the solution of problems. Not necessarily the smartest and shortest path. Not immediately, anyhow.

A possible different teaching sequence would be to start from problems: give the students some problems, which once mathematized turn out into solving a second degree polynomial. Some of the problems should be aesily solvable (i.e. lead to an equation of the form $x^2=d$ or $x^2+bx=0$), some should lead to a complete equation with no vanishing terms. The students will find by themselves how to solve the simpler ones, and maybe will even give a try to the more difficult ones. Guided by the teacher, working in groups, they may rediscover themselves the formula (maybe by completing the square) and teach it to other groups. A discover made by themselves, while trying to solve a problem and effectively experiencing the hardness of the problem and the sense of joy that comes with the solution, will leave the students something more than just a formula to blindly apply. Most of all, it will leave them with the sense of doing math.

And after the group work, a good recap by the teacher would be nice, so to put all the ideas (which came from the students) in order. Doing like that, probably a lot more of them will remember the formula, but most important will know how to find it again if needed.

\begin{question}{So... is math usefulness relevant to learners?} In my opinion it is very little relevant for their decision to willingly learn maths, and both teachers and popularizers of mathematics should not focus too much on applications and usefulness of mathematics, since applications and usefulness are often not immediate, but rather on the joy and the challenges of doing mathematics.
\end{question}

\section{Is math popularization useful? --- A math popularizer's apology}
To end this chapter, I would like to give an apology to the activity of math popularization. Hardy had quite harsh words for this activity, and I feel most of my collegues agree with him: any amount of time spent in popularizing math is time stolen from doing actual math, and probably, if you do that, it is just because you are not good enough to do actual math.

\begin{quotation} If then I find myself writing, not mathematics, but ‘about’
mathematics, it is a confession of weakness, for which I may
rightly be scorned or pitied by younger and more vigorous
mathematicians. I write about mathematics because, like any
other mathematician who has passed sixty, I have no longer the
freshness of mind, the energy, or the patience to carry on
effectively with my proper job.
\begin{flushright}A Mathematician's Apology \cite{H} \S 1 --- G. H. Hardy
\end{flushright}
\end{quotation}

I'm sorry for my colleagues, but if this vision of popularization or communication of mathematics was ok in the Fourties of the last century, it is no longer so nowadays. For a better and longer essay on this subject, I refer the reader to the article by Silvia Benvenuti and Roberto Natalini \cite{BN}.

There are several top mathematicians deeply involved in communication of math (just think at the Fields medallists Cedric Villani and Alessio Figalli, to name two). Moreover, in the 2011 \textit{European's Charter of Researchers} \cite{ECR} it is clearly written that scientists should be directly involved in communicating their own researches to the wide public in order to favour the creation of a scientific mind.

\begin{quotation}Researchers should ensure that their research activities are made known to society at large in such a way that they can be understood by non-specialists, thereby improving the public’s understanding of science. Direct engagement with the public will help researchers to better understand public interest in priorities for science and technology and also the public’s concerns.
\begin{flushright}European's Charter of Researchers \cite{ECR}
\end{flushright}
\end{quotation}

The reason for that is precisely what we tried to outline and suggest in this chapter: a scientific-leaned mind is needed for the well-being of the society at large, and ---given how modern democracies work--- it is a need of the whole society and not just for a few enlightened who are part of the governing class.

The idea many researchers have of themselves and of research is that they are needed by society (which is true) and they have no urge to explain to society why they are needed (and this is false). This idea that what matters for research is getting it done and not being presented to society at large is deeply fixed into researchers' minds, but it is false, in the sense that the society must be aware of the fact that investing (money, time and effort) in research, both applied and pure, is what we need to do. And this is even more true when we talk of an inherently abstract subject as mathematics, whose practical implications are neither immediate nor evident.

\begin{figure}[h]
\center
		\includegraphics[width=0.55\textwidth]{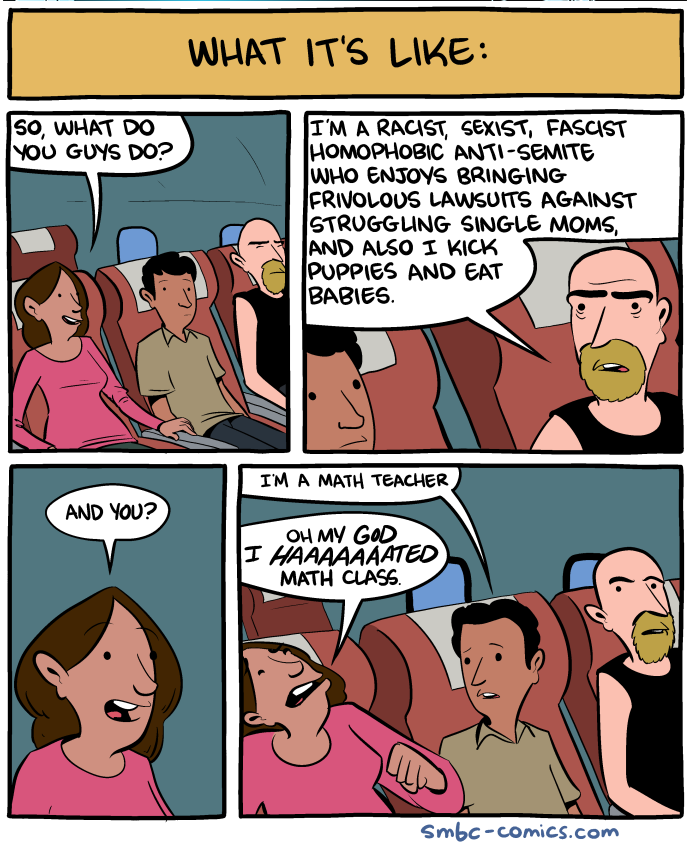}
	\caption{A comic by Saturday Morning Breakfast Cereal \cite{SMBC} getting right to the point of how math is perceived.}
	\label{SMBC}
\end{figure}

I perfectly know the feeling of frustration when you are telling someone you are a mathematician, or you teach math and the response you get is that shown in the comic by SMBC (see figure \ref{SMBC}).

You usually get on the defensive and have trouble to communicate the beauty of mathematics, or even ---if you are tired--- do not gey at all into the subject. Or sometimes, people just say that they understand math is useful, but \textit{not for them} / \textit{they do not get it} / \textit{they are not a math person} (chose one or more).

Communicating to the society is tough. And society at large are not people who willingly go to an event of science (or math) popularization, or not just them: society at large, like it or not, is mainly composed by people who have a problematic relationship with mathematics, and they will not come to an event where you talk math to them.

\begin{quotation} Part of the problem is that nobody has the faintest idea what it is that mathematicians do.
\begin{flushright}A Mathematician's Lament \cite{L}  --- P. Lockhart
\end{flushright}
\end{quotation}

Doing mathematics is an actvity quite similar to that of an artist or a writer: there is a lot of technique involved, but also a lot of artistic out-of-the-box thinking and imagination. People are scared by technique (the only thing about math they know) and are not willing to know more about math.

\begin{figure}[h]
\center
		\includegraphics[width=0.70\textwidth]{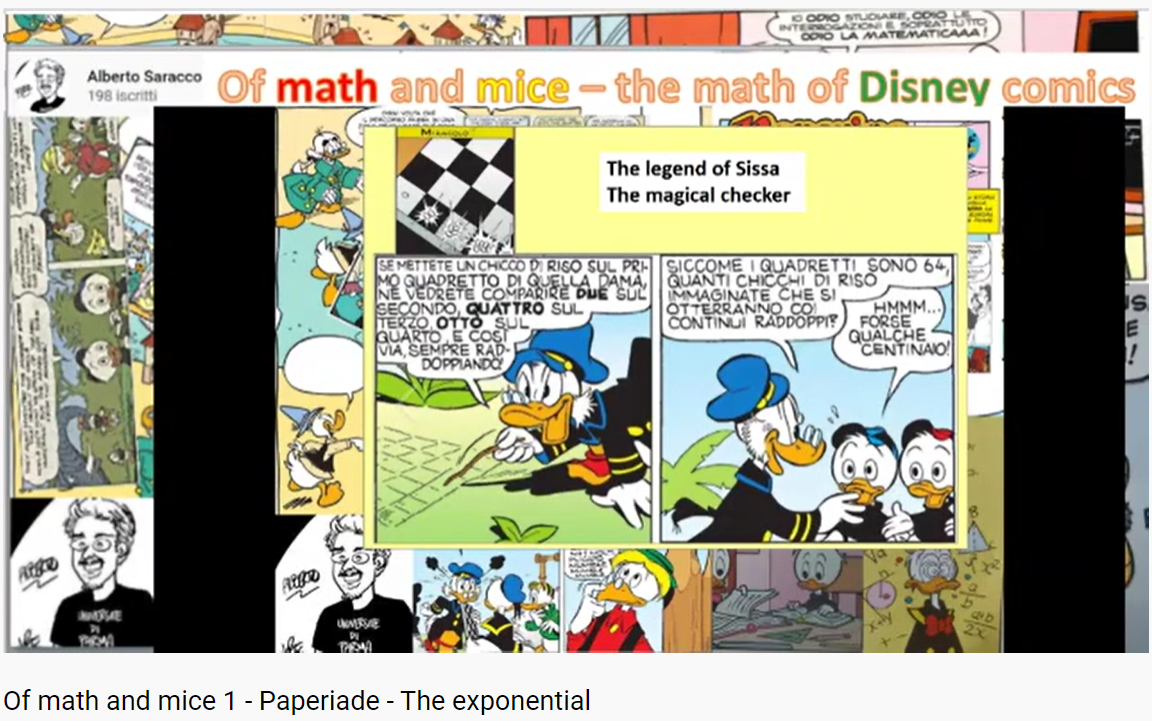}
	\caption{The first video of the YouTube project \textit{Of math and mice - the mathematics of Disney comics} \cite{OMAM}, devoted to popularizing math using Disney comics, translation in English of the corresponding Italian project \textit{Un matematico prestato alla Disney} \cite{UMPAD}, both available on my YouTube channel \textit{Alberto Saracco}.}
	\label{OMAM1}
\end{figure}

It is up to you as a mathematician to get people to know what mathematicians do. They won't come at you. You have to come at them, by using their passions to talk about math. It is a while, since I started doing that with comics, using Disney comics to talk about math (see \cite{S}, but also my YouTube playlists on the subject \cite{OMAM,UMPAD}, see figure \ref{OMAM1}).

I noticed that when I put Disney comics or characters in the title of one of my talks, the audience is thrice as big as it usually is. And often many of the people in the audience were not interested at all about math at the beginning of the talk, but exited the talk with better feelings about the subject.

We need these kind of things in order to get math near to people who would not approach math. And, as I said, it is of fundamental importance. If we want society to grow and fully use math, it is up to us. As Francesca Arici, member of the \textit{Raising Public Awareness} European Commettee, said

\begin{quotation} Don't be afraid to use metaphors and don't be afraid to lie a little bit: the objective is to communicate math and not to make people see we can do math and prove theorems in a rigorous way.
\begin{flushright} F. Arici \cite{Ar}
\end{flushright}
\end{quotation}


\begin{thebibliography}{99.}%
%
%
\bibitem{A}Abate, M.: Scrivere Matematica nel fumetto, In: Emmer, M. (ed.)  Matematica e cultura 2004, pp. 19–29, Springer Italia, Milano (2004)

\bibitem{AG} Abstruse Goose: \textit{Impure Mathematics}, 
abstrusegoose.com/504.
\bibitem{Ar} Arici, F. and Mulas, R.: \textit{Matematica, Teatro e Metafore, con Francesca Arici - La matematica danzante}, YouTube video on the channel \textit{Madd Maths} (2021)

\bibitem{KA} Arrow, K. J., \textit{Social Choice and Individual Values}, Cowles Commission Monograph No. 12, John Wiley \& Sons, Inc., New York, 1951.

\bibitem{B} Bazett, T.: \textit{How SIMPSON'S PARADOX explains weird COVID19 statistics}, YouTube video on the channel \textit{Dr. Trefor Bazett} (2020)

\bibitem{BN} Benvenuti, S. and Natalini, R.: Comunicare la matematica: chi, come, dove, quando e,
soprattutto, perchè?! Mat. Cult. Soc. \textbf{2}-2, 175--193 (2017)

\bibitem{H} Hardy, G. H., \textit{A Mathematician's Apology}, Cambridge University Press, Cambridge, 1940.\\ Available on-line: www.math.ualberta.ca/mss/misc/A\%20Mathematician\%27s\%20Apology.pdf 

\bibitem{L} Lockhart, P., \textit{A Mathematician's Lament}, Bellevue Literary Press, 2009. The first version (2002) is available on-line: mysite.science.uottawa.ca/mnewman/LockhartsLament.pdf

\bibitem{MP} Minute Physics: \textit{Simpson's Paradox}, YouTube video on the channel \textit{minutephysics} (2017)

\bibitem{NG} Nicolai, R., Gottardo, A.: Paperone, Rockerduck e i consigli incrociati. Topolino \textbf{2945}-2, 25pp. (2012)

\bibitem{GB} Poe, E. A.: \textit{The Gold-Bug}, Dollar Newspaper (Philadelphia, PA), vol. I, no. 23, June 28, 1843, pp. 1 and 4. Available online: www.eapoe.org/works/tales/goldbga2.htm

\bibitem{S}Saracco, A.: Math in Disney Comics, In: Abate, M.,  Emmer, M. (eds.) Imagine Math 7, pp. 189–210, Springer Italia, Milano (2020)

\bibitem{OMAM} Saracco, A.: \textit{Of math and mice - the mathematics of Disney comics}, YouTube playlist on the channel \textit{Alberto Saracco} (2020) 

\bibitem{UMPAD} Saracco, A.: \textit{Un matematico prestato alla Disney}, YouTube playlist on the channel \textit{Alberto Saracco} (2020)

\bibitem{SS19} Saracco, A. and Saracco, G.: A discrete districting plan. Netw. and Heterog. Media, \textbf{14}(4), 771--788 (2019)

\bibitem{SS} Saracco, A. and Saracco, G.: Matematica ed elezioni, paradossi e problemi elettorali. Mat. Cult. Soc. \textbf{5}-1, 17--31 (2020)

\bibitem{SMBC} Saturday Morning Breakfast Cereal: \textit{What it's like}, 
www.smbc-comics.com/comic/what-its-like. 

\bibitem{RV} Sisti, A., Intini, S.: 3 gradini per le stelle - Paperino Paperotto, Roby Vic e i conti... alla rovescia. Topolino \textbf{2859}-1, 34pp. (2010)

\bibitem{Z} Zaccagnini, A.: Dialogo sui numeri primi - un dialogo galileiano (2020) e-book available on-line:
maddmaths.simai.eu/divulgazione/zaccagnini-dialogo-ebook/

\bibitem{ECR} \textit{European's Charter of Researchers}, 
euraxess.ec.europa.eu/jobs/charter/european-charter.

\end{thebibliography}
\end{document}